\documentstyle[11pt,amssymb]{article}

\begin{document}

\hsize=15cm
\vsize=20cm

\newcommand{\f}{ \cal F}
\newcommand{\h}{{\sf H}}

\newtheorem{De}{Definition}[section]
\newtheorem{Th}[De]{Theorem}
\newtheorem{Pro}[De]{Proposition}
\newtheorem{Le}[De]{Lemma}
\newtheorem{Ex}[De]{Example}
\newtheorem{Co}[De]{Corollary}

\newtheorem{dummy}[De]{}

\newcommand{\Def}[1]{\begin{De}#1\end{De}}
\newcommand{\Thm}[1]{\begin{Th}#1\end{Th}}
\newcommand{\Prop}[1]{\begin{Pro}#1\end{Pro}}
\newcommand{\Exa}[1]{\begin{Ex}#1\end{Ex}}
\newcommand{\Lem}[1]{\begin{Le}#1\end{Le}}
\newcommand{\Cor}[1]{\begin{Co}#1\end{Co}}

\centerline {\bf  VANISHING LINE FOR THE DESCENT SPECTRAL SEQUENCE }

\bigskip
\bigskip

\centerline{\sc {Teimuraz Pirashvili}}
\bigskip
\bigskip
\centerline{\it A.M. Razmadze Math. Inst., Aleksidze str. 1, }
\centerline {\it Tbilisi, 380093. Republic of Georgia}
\centerline {{\it email pira@mpim-bonn.mpg.de} and {\it pira@rmi.acnet.ge}}
\bigskip
\bigskip

\noindent {\bf Abstract}.  {\it We prove that if $f:X\to Y$  is a 
surjective and proper map between
locally compact spaces and if additionally it is  
cohomologically $(k-1)$-connected, 
for some $k>0$, then for the cohomological descent 
spectral sequence \cite{de} one has $E^{pq}_2=0$ provided $q<pk$.}



\bigskip
\bigskip
\noindent {\bf 1. Introduction}. Let $f:X\to Y$ be a 
continuous map of topological spaces and let 
${\cal F}$ be an abelian sheaf on $Y$. 
One defines the simplicial space $X_*$ by
$$X_n=\lbrace (x_0,...,x_n)\in X^{n+1}\mid f(x_0)=...=f(x_n)\rbrace ,n\geq 0, $$
$$\partial _i(x_0,...,x_n)=
(x_0,...,{\widehat x_i},...,x_n),0\leq i\leq n$$
$$s_i(x_0,...,x_n)=(x_0,...,x_i,x_i,...,x_n),0\leq i\leq n.$$
We let $f_n:X_n\to Y$ be the map given by $f_n(x_0,...,x_n)=f(x_0)=...=f(x_n)$.
Moreover we put ${\f}_n=f_n^*{\f}$. 
Thus $X_0=X$,  $f_0=f$ and ${\f}_0=f^*{\f}$. Then $(f_n)_{n\geq 0}$ 
defines the morphism 
from the simplicial space $X_*$ to the constant simplicial space 
$Y$ and $\f _*$ is a sheaf on the simplicial space $X_*$. 
The following spectral 
sequence is known as the cohomological descent spectral 
sequence (see for example
5.3 of \cite{de} and references given there).

\bigskip
\noindent {\bf Theorem 1.} {\it Let $f:X\to Y$ be a 
surjective and proper map, then there exists a 
spectral sequence
$$E_1^{pq}=\h^q(X_p,{\cal F}_p)\Longrightarrow \h^{p+q}(Y,{\cal F}).$$
Moreover $E^{pq}_2$ is isomorphic to the 
$p$-dimensional cohomology of the cochain complex associated 
to the cosimplicial abelian group  $E_1^{*q}$.
}

The aim of the paper is to prove the following vanishing result for the 
cohomological descant spectral sequence. In this paper all maps between
topological spaces are continuous and all sheafs are abelian. Moreover a 
space $Z$ is 
called {\it cohomologically $k$-connected} if for any constant 
sheaf $A$ one has
$$\h^0(Z,A)\cong A, \ \ {\rm and} \ \ 
\h^i(Z,A)=0 \ \ {\rm for }\ \ 0<i\leq k.$$
Similarly, a map $f:X\to Y$ is called {\it cohomologically $k$-connected} if for each $y\in Y$ the space $f^{-1}(y)$ is cohomologically $k$-connected.

\bigskip
\noindent {\bf Theorem 2.} {\it Let $f:X\to Y$  be a 
surjective and proper map between
locally compact spaces. If $f$ is  cohomologically $(k-1)$-connected, 
for some $k>0$, then for the cohomological descent 
spectral sequence one has $E^{pq}_2=0$ provided $q<pk$.}

\smallskip
\noindent Actually we will prove much strong fact: the normalization of the 
cosimplicial abelian group $E^{*q}_1$ vanishes in dimension $p$ 
as soon as $q<kp$.

\smallskip
\noindent We refer the reader to \cite{P} for different sort
of vanishing results in group cohomology framework.

\bigskip
\noindent The author was  partially supported by the grant
INTAS-99-00817  and by the TMR network K-theory and algebraic groups, 
ERB FMRX CT-97-0107. Special thanks to MPI at Bonn for hospitality.

\bigskip
\noindent {\bf 3.  Preliminaries on cosimplicial abelian groups}.
Let  $(A^*,\delta^*,s^i)_{0\leq i\leq n}$ be a cosimplicial abelian group. 
It can be considered as 
a cochain complex with the coboundary map given by 
$d=\sum_i(-1)^i\partial ^i$. The normalization of $A^*$ is 
a subcomplex $N^*A^*$, given by
$$N^n(A^*)=\{ x\in A^n\mid s^ix=0,i=0,\cdots ,n-1\}.$$
It is well-known that $N_nA^*$ is isomorphic to the quotient of 
$$A_n /(Im (\delta ^0)+ \cdots + Im (\delta ^{n-1}))$$ and 
the inclusion $N^*A^*\subset A^*$ is a 
homotopy equivalence. Thus $H^*(A^*,d)\cong H^*(N^*A^*)$.
Moreover a theorem of Dold and Kan asserts that the functor 
$N^*$ establishes an equivalence between the category of cosimplicial 
abelian groups and the category of nonnegative cochain complexes. In particular
$N^*$ is an exact functor and therefore 
for any short exact sequence of 
cosimplicial abelian groups
$$0\longrightarrow A^*_1\longrightarrow A^*
\longrightarrow A_2^*\longrightarrow 0$$
one has
$$N^i(A^*)= 0,\ \ {\rm iff }\ \ N^i(A^*_{1})=0=N^i(A_2^*) .$$
Now we consider the following example. Let $R$ be an algebra over
a commutative ring $K$. We let $C^*(R)$ be the following 
cosimplicial $K$-module 
$$C^n(R)=R^{\otimes (n+1)},$$
$$\delta ^i(x_0,\cdots ,x_n)=(x_0,\cdots, x_{i-1},1,x_i,\cdots,x_n),$$
$$s^i(x_0,\cdots,x_n)=(x_0,\cdots, x_ix_{i+1},\cdots, x_n).$$
It is clear that for any $R$ one has
$$N^n(C^*(R))={\bar R}^{\otimes n}\otimes R,$$
where ${\bar R}=R/K.1.$ If $R$ is graded $R=\bigoplus_{i\geq 0}R^i$, then 
one has an isomorphism of  cosimplicial modules 
$C^*(R)\cong \bigoplus_{m\geq 0}C^{*m}$, where $C^{nm}= 
\bigoplus R^{ i_0}\otimes \cdots \otimes R^{i_n}$
and the  sum is taken over all 
$i_0+\cdots +i_n=m, i_s\geq 0.$ It is now clear that, 
if $R$ is graded and 
connected (that is $R^0=K$), then $N^nC^{*m}=0$ if $m<n$. More generally, 
we have $N^nC^{*m}=0$ if $m<kn$, as soon as $R$ is $(k-1)$-connected 
for some $k\geq 1$ (that is connected and $R^i=0$ for $0<i<k$). 

We will use this observation in the following situation. Let $S$ be a compact
space, which is assumed to be cohomologically $(k-1)$-connected. We let 
$c_*(S)$ be the simplicial space given by
$$c_n(S)=S^{n+1},$$ 
$$\partial _i(x_0,\cdots ,x_n)=(x_0,\cdots ,\hat x_i,\cdots, x_n),$$
$$s_i(x_0,\cdots ,x_{n-1})=(x_0,\cdots ,x_i,x_i,\cdots ,x_{n-1}).$$
Now applying the functor $\h ^m(-,A)$ we get a cosimplicial 
abelian group $\h^m(c_*(S),A)$.
Here $A$ is an abelian group, considered as
a constant sheaf.  We claim that the normalization of this 
cosimplicial abelian group vanishes in dimension $n$ provided $m<nk$. 
 It is well known that  the cohomology of a compact space with  
coefficients in a constant sheaf is isomorphic 
to the  {\v C}ech cohomology and 
therefore one has the following natural short 
exact sequence
$$0\to \h^i(S,{\bf Z})\otimes A\to \h^i(S,A)\to \h^{i-1}(S,{\bf Z})*A\to 0. \leqno (4)$$
Indeed, since the tensor and the torsion products preserves colimits, 
it suffice to
consider the case, when $S$ is a finite polyhedron. In this case 
the {\v C}ech cohomology is isomorphic 
to the simplicial cohomology, and hence $\h^*(S,A)$ 
is the homology of ${\sf Hom}(C_*(S),A)$ with 
degreewise free and finitely generated $C_*(S)$. 
Now it is enough to observe that the natural map
$$ {\sf Hom}(C_*(S),{\bf Z})\otimes A\to {\sf Hom}(C_*(S),A)$$
is an isomorphism. It follows from the exact sequence (4) 
that to prove the claim it suffice to consider the case, when
$A=K$ is a field. In this case one can use  
the K\"unnent theorem for ${\rm {\hat C}}$ech cohomology 
(see for example Theorem 7.2 of \cite{mas}) to obtain the isomorphism 
that $C^*(\h^*(S,K))\cong \h^*(c_*(S),K).$ 
Thus $N^m(\h^n(c_*(S),K)=0$
provided $n<mk$ and the claim is proved.

\bigskip
\noindent{\bf 5. Preliminaries on the Lerey spectral sequence}.
 Let $g: W\to Z$ be a map of topological spaces and let 
${\cal G}$ be a sheaf on $W$. Then there exists a spectral sequence
$$E^{pq}_2=\h^p(Z, {\cal R}^qg_*({\cal G}))\Longrightarrow 
\h^{p+q}(W,{\cal G})$$
known as the Lerey spectral sequence. Here ${\cal R}^qg_*(\cal G)$ is 
the sheaf on $Z$ 
associated to the presheaf $U\mapsto \h^q(g^{-1}(U),{\cal G})$. Here $U$ is 
an open subset of $Z$. It is well-known (see page 201 of \cite{gode}) that 
$${\cal R}^q{g_*(\cal G)}_z\cong \h^q(g^{-1}(z),{\cal G})$$
as soon as $g$ is proper and $W,Z$ are locally compact. 
Now we take ${\cal G}=g^*{\cal B}$, where ${\cal B}$ is a sheaf on $Z$. 
Then for each $w\in W$ one has ${\cal G}_w={\cal B}_{g(w)}$. 
Thus for each $z\in Z$ the sheaf $\cal G$ is a constant sheaf on $g^{-1}(z)$ associated to the abelian group ${\cal B}_z$. In this case we
 allow ourself to rewrite the Leray spectral sequence as follows
$$E^{pq}_2=\h^p(Z, z\mapsto \h^q( g^{-1}(z), {\cal B}_z))\Longrightarrow 
\h^{p+q}(W,g^*{\cal B}).$$

\bigskip
\noindent{\bf 6. Proof of Theorem 2}. For each $n\geq 0$ we 
consider the Leray 
spectral sequence corresponding to $f_n$:
$$_nE^{pq}_2=\h^p(Y, y\mapsto \h^q( f_n^{-1}(y), {\cal F}_y)\Longrightarrow 
\h^{p+q}(X_n, {\cal F}_n).$$
Varying $n$ one obtains the spectral sequence of cosimplicial abelian groups.
Since the normalization is an exact functor it suffices to show that for 
fixed $p$ and $q$ the normalization of the cosimplicial abelian group
$[n]\mapsto\  _nE^{pq}_2$ vanishes in the dimension $m$, provided $q<mk$. 
Let us consider the cosimplicial object in the category of sheafs over $Y$ given by
$$[n]\mapsto {\cal R}_*^q(f_n)_*(f_n^*\f),$$
which we denote for simplicity by $[n]\mapsto ({\cal R}^q)^n$. Thus 
$$_*E^{pq}_2=\h^p(Y, ({\cal R}^q)^*).$$
It follows from the exactness of the normalization that
$$N^m(_*E^{pq}_2)=\h^p(Y,N^m({\cal R}^q)^*).$$
Therefore it suffices to show that $N^m({\cal R}^q)^*=0$ provided $m<qk$. 
Since the stalk at $y$ is an exact functor, we have 
$$(N^m({\cal R}^q)^*)_y\cong N^m(({\cal R}^q_y)^*).$$
The last group is the same as $N^m\h^q(f_*^{-1}(y),{\f}_y)$. 
It suffice to show that it vanishes as 
soon as $q<km$. But this follows from the fact that
$$f_n^{-1}(y)=f^{-1}(y)\times \cdots \times f^{-1}(y)$$
(and therefore $f_*^{-1}(y)=c_*(f^{-1}(y))$) and from the 
example of Section 3.

\end{document}